\pdfoutput=1
\documentclass[10pt,twocolumn,twoside]{IEEEtran}
\usepackage{pdfpages}
\usepackage{graphicx}
\usepackage{epsfig}
\usepackage{amssymb}
\usepackage{enumerate}
\usepackage{latexsym}
\usepackage[cmex10]{amsmath}
\newtheorem{definition}{Definition}

\newtheorem{remark}{Remark}
\newtheorem{theorem}{Theorem}

\def\pdH{{ {\cal H}^+}(p)}

\def\bU{{\boldsymbol U}}

\def\diag{ {\rm diag} }

\def\EqualDist{\mathrel{\mathop=^{\rm d}}}
\def\EqualDef{\mathrel{\mathop=^{\rm def}}}

\def\bSigma{{\boldsymbol\Sigma}}

\def\bLambda{{\boldsymbol\Lambda}}

\def\bX{{\boldsymbol X}}

\def\bQ{{\boldsymbol Q}}
\def\bA{{\boldsymbol  A}}

\def\bC{{\boldsymbol  C}}

\def\bB{{\boldsymbol  B}}

\def\bI{{\boldsymbol  I}}

\def\bL{{\boldsymbol  L}}

\def\bP{{\boldsymbol  P}}

\def\bR{{\boldsymbol  R}}

\def\bS{{\boldsymbol  {S}}}

\def\tr{{\rm tr}}

\def\half{{\textstyle{\frac 1 2}}}

\def\a{ {\rm a}}

\def\EqualDist{\mathrel{\mathop=^{\rm d}}}

\def\sxy*{\mathbf{S}_{xy^{*}}(f)}
\def\syx*{\mathbf{S}_{y^{*}x}(f)}

\def\var{ {\rm var} }

\begin{document}

\title{Intrinsic Analysis of the Sample Fr{\'e}chet Mean and Sample Mean of Complex Wishart Matrices}

\author{L.~Zhuang \&  A.~T.~Walden,~\IEEEmembership{Senior~Member,~IEEE}
\thanks{
Copyright (c) 2017 IEEE. Personal use of this material is permitted. However, permission to use this material for
any other purposes must be obtained from the IEEE by sending a request to pubs-permissions@ieee.org. 
Linjie Zhuang and Andrew Walden
are both at the Department of Mathematics, Imperial College  London, 180 Queen's Gate,
London SW7 2BZ, UK.  
(e-mail:  l.zhuang13@imperial.ac.uk and a.walden@imperial.ac.uk).} }
\IEEEpubid{}
\maketitle
\begin{abstract}
We consider two types of averaging of complex covariance matrices, a  sample mean (average) and the sample Fr{\'e}chet  mean. We analyse the performance of these quantities as estimators for the true covariance matrix via `intrinsic' versions of bias and mean square error, a methodology which takes account of geometric structure. We derive simple expressions for the intrinsic bias in both cases, and the simple average is seen to be preferable. The same is true for the asymptotic Riemannian risk, and for the  Riemannian risk itself in the scalar case. Combined with a similar preference for the simple average using non-intrinsic analysis, we conclude that the simple average is preferred overall to the sample Fr{\'e}chet  mean in this context.
\end{abstract}

\section{Introduction}
In this letter we give a performance analysis of the average of complex-valued covariance matrices found by the simple sample mean or via the sample Fr{\'e}chet mean. Classical measures such as bias and mean square error (MSE) do not take into account the geometrical structure.

 Covariance
matrix averaging methods for EEG signal classification was of interest in 
\cite{Yger_etal15}, and robust averaging was proposed in \cite{Ilea_etal16} for the case of inhomogeneous samples due to outliers. Our motivation to study this topic came from the desire to build improved graphical models from neuroscience data, 
used in schizophrenia studies \cite{Medkouretal10}; see also \cite{Medkouretal09}. 
This requires good spectral matrix estimates for different groups of individuals and we want to know whether it is best
to average estimated spectral matrices over individuals in each group, or find instead their sample
Fr{\'e}chet mean.

Positive definite and Hermitian complex
covariance matrices  form a manifold in the space of complex-valued matrices.
If such a manifold is equipped with a Riemannian metric, it becomes a Riemannian manifold. In evaluating performance criteria for the different forms of averaging of covariance matrices, the manifold should be taken into account. This leads to considering `intrinsic' versions of bias and MSE \cite{GarciaOller06,OllerCorcuera95,Smith05}.
Indeed this was the approach adopted for performance analysis by
\cite{Ilea_etal16}.

In a recent work \cite{ZhuangWalden17} the averaging of 
a set $\bS_1,\ldots,\bS_N$ of positive definite complex covariance matrix estimators  was studied. The matrices  are considered homogeneous in the sense that they are  estimating the same true covariance matrix $\bSigma,$ and statistically, they are independent and identically complex-Wishart distributed. Ordinary bias and MSE were used to compare the performance of the ordinary 
average of the estimated covariance matrices, and of their  sample Fr{\'e}chet mean, (also called the Riemannian mean in this setting). In this letter 
we derive the intrinsic versions of bias and MSE to enable a
`geometry-aware' performance appraisal.
In contrast to the classical measures, intrinsic evaluation of an estimator is invariant under reparametrization, and it is dependent on the model only
\cite{GarciaOller06}.

The main contribution of this letter is to prove
the simple average of $N$ Wishart matrices outperforms their Fr\'echet mean 
in terms of intrinsic bias for $N\geq 2.$ We also show that the same is true for the asymptotic Riemannian risk, and for the  Riemannian risk itself in the scalar case for $N\geq 2.$ These results are enabled by the derivation of simple expressions for the intrinsic bias for both types of average.

\section{Preliminaries}\label{sec:prelim}

\subsection{Complex Wishart matrices}\label{subsec:CWish}
Let ${\bX}_0,\ldots,{\bX}_{K-1}$ be $K$ independent $p$-dimensional complex-Gaussian  random vectors  with zero means and covariance matrix $\bSigma.$
Then the maximum likelihood estimator for $\bSigma$ is the covariance matrix estimator
$
{\bS}=(1/K) \sum_{k=0}^{K-1} \bX_k\bX_k^H,
$
where $^H$ denotes Hermitian transpose
and $K{\bS}$ has the complex Wishart distribution \cite{Goodman63} with $K$ complex degrees of freedom and mean $K\bSigma,$ denoted by
\begin{equation}\label{eq:Wishart}
K {{\bS}}
 \EqualDist
{\cal W}_p^C(K,{\bSigma}).
\end{equation}
Such matrices arise frequently, e.g., \cite{Conradsen_etal03,Medkouretal10}. We assume $K\geq p,$ since then $\bS$ has full rank $p,$ and $\bS$ is positive definite, with probability one.

\subsection{Positive Definite Hermitian matrices}\label{sec:posdef}
Let $\pdH$ denote the set of $p \times p$ positive definite Hermitian matrices,  a differentiable manifold.  
Denote the set of all $p \times p$ invertible matrices by $GL(p)$. The group action of $GL(p)$ on  $\pdH$ is the transformation
$\phi: GL(p) \times \pdH \rightarrow \pdH$ given by $\phi_\bL(\bP) = \bL \bP\bL^H$. 

\subsection{Riemannian manifold}
 A manifold equipped with a Riemannian metric $g$ is a Riemannian manifold. 
$\pdH$ can be turned into a Riemannian manifold by defining at every point in $\pdH$ an inner product for elements in the tangent space that varies differentiably along the manifold.
 The minimum length curve connecting two points on the manifold is called the geodesic and the Riemannian distance $d_g$ between the points is given by the length of this curve. We  use a scaled version of the Frobenius inner product that is invariant under the group action of $GL(p).$ 
For any $\bA,\bB\in T_{\bP}\pdH$, 
where $T_{\bP}\pdH$ denotes the tangent space at the point $\bP\in \pdH,$
the `scaled Frobenius inner product' takes the form (e.g., \cite[p.~1729]{Jiang_etal12b})
\begin{equation}\label{eq:scFrob}
\left\langle\bA,\bB\right\rangle_{\bP} =\text{tr}(\bP^{-1}\bA\bP^{-1}\bB),
\end{equation}
and the corresponding norm is 
\begin{equation}\label{eq:wantc}
\|\bA\|_{\bP}=\left\langle\bA,\bA\right\rangle _{\bP}^{1/2}=||\bP^{-1/2}\bA\bP^{-1/2}||_{\rm Fr}.
\end{equation}

\subsection{Riemannian distance}
For $\bP_0 \in GL(p)$ the matrix logarithm is any $p \times p$ matrix $\bQ$ such that $\exp(\bQ)=\bP_0,$ where $\exp(\cdot)$ is the usual power series expansion.
For $\bP_1\in \pdH$ having spectral decomposition $\bP_1=\bU \bLambda \bU^H,$ 
the matrix logarithm function is 
 \cite[p.~429]{Bernstein05}
\begin{eqnarray}
\!\!\!\!\!\log\bP_1 
&=&\bU\diag(\log\lambda_1(\bP_1),\ldots,\log\lambda_p(\bP_1))\bU^H,
\label{eq:deflogm}
\end{eqnarray}
where $\lambda_j(\bP_1)>0$ is the $j$th eigenvalue of $\bP_1.$ 
From (\ref{eq:deflogm})  $(\log\bP_1)^H=\log\bP_1,$ so  $\log \bP_1$ is Hermitian. For any invertible matrix $\bA$ and  matrix $\bB$ having real positive eigenvalues, 
\begin{equation}\label{eq:onec}
\log(\bA^{-1}\bB\bA)=\bA^{-1}\log(\bB)\bA.
\end{equation} 
The Riemannian distance $d_g(\bP_1,\bP_2)$ is defined by \cite{Jiang_etal12b} 
\begin{eqnarray}
d_g^2(\bP_1,\bP_2) &=& \tr(\log^2(\bP_1^{-1/2}\bP_2\bP_1^{-1/2}))\label{eq:cdista}\\
&=& \|\log(\bP_1^{-1/2}\bP_2\bP_1^{-1/2})\|^2_{\rm Fr}.
\label{eq:cdistb}
\end{eqnarray}

We will denote the Riemannian manifold by ${\cal M}.$

\subsection{Fr{\'e}chet  Mean}
Let $\bS_1,\ldots,\bS_N$ be independent random matrices with common distribution $F$ on  ${\cal M}.$ Let ${\hat F}$ be the empirical distribution.
The sample
Fr{\'e}chet  mean  of  ${\hat F}$ is the  minimizer of \cite{BhattacharyaPatrangenaru03}
\begin{equation}\label{eq:posdef}
(1/N)\sum_{j=1}^N d_g^2({\hat\bP},\bS_j), \quad {\hat\bP}\in \pdH.
\end{equation}
The Riemannian metric space has negative sectional curvature so the sample Fr{\'e}chet mean is unique, \cite[p.~6]{Dryden_etal09}.

\subsection{Mappings for Riemannian Manifolds}
The exponential map $\text{Emap}:T_{\bP}{\cal M}\longrightarrow {\cal M}$ is a function mapping a vector $\bU$ (starting from $\bP\in {\cal M}$) in the tangent space, to a point $\bS$ on the Riemannian manifold:
\begin{equation*}
\bS=\text{Emap}_{\bP}(\bU)=\bP^{1/2}\exp(\bP^{-1/2}\bU\bP^{-1/2})\bP^{1/2}.
\end{equation*}
 The (inverse) logarithmic map is:
\begin{equation*}
\bU=\text{Lmap}_{\bP}(\bS)=\bP^{1/2}\log(\bP^{-1/2}\bS\bP^{-1/2})\bP^{1/2},
\end{equation*}
and takes $\bS$ on the Riemannian manifold to $\bU$ in the tangent space, $
\mathcal{M}\longrightarrow T_{\bP}\mathcal{M}
$
;
the nonpositive curvature  of the manifold guarantees the unique inverse mapping.

\section{Intrinsic Versions of Bias and MSE}

Before proceeding we need to state some definitions.

\begin{definition}\cite[p.~128]{GarciaOller06}
For a fixed sample size $N,$ and an estimator
${\cal S}$ of $\bSigma,$ the estimator vector field is 
$\text{Lmap}_{\bSigma}({\cal S}).$
\end{definition}
\begin{definition}
\cite[p.~1615]{Smith05}
The expectation of the estimator ${\cal S}$ on the manifold,
with respect to $\bSigma,$  denoted here
${\cal E}_\bSigma\{ {\cal S} \}\in {\cal M},$ is defined as
$
{\cal E}_\bSigma\{ {\cal S} \}=\text{Emap}_{\bSigma}(E\{ \text{Lmap}_{\bSigma}({\cal S})\}).
$
\end{definition}
\begin{remark}
An estimator ${\cal S}$ is an {\it intrinsically unbiased} estimator of $\bSigma$ if and only if $E\{ \text{Lmap}_{\bSigma}({\cal S})\}={\bf 0},$ since then
$ {\cal E}_\bSigma\{ {\cal S} \}=\bSigma^{1/2}\exp({\bf 0}) \bSigma^{1/2}=\bSigma.$ In fact $E\{ \text{Lmap}_{\bSigma}({\cal S})\}$ is called the bias vector field, e.g., \cite[p.~1615]{Smith05}.\hfill$\square$
\end{remark}
\begin{definition}\label{def:int_bias}
\cite[p.~129]{GarciaOller06}, \cite[p.~1568]{OllerCorcuera95}
The intrinsic (or invariant) bias of an estimator ${\cal S}$ of $\bSigma$ is a scalar defined as
\begin{equation}\label{eq:defibias}
\text{ibias}({\cal S})\EqualDef \|E\{ \text{Lmap}_{\bSigma}({\cal S})\}\|^2_{\bSigma}.
\end{equation}
\end{definition}
Hence, calculating the intrinsic bias of an estimator involves (i) mapping the estimator to a vector on the tangent space specified by the true parameter; (ii) finding the squared norm of the expectation of the vector. 
\begin{definition}
\cite[p.~129]{GarciaOller06}, \cite[p.~1568]{OllerCorcuera95}
The Riemannian risk of ${\cal S}$ is a scalar defined as 
\begin{equation}\label{eq:defRisk}
E\{ \|\text{Lmap}_{\bSigma}({\cal S})\|^2_{\bSigma}\}.
\end{equation}
\end{definition}
\begin{remark}
The Frobenius inner product is a Riemannian metric in the Euclidean domain and makes $\mathcal{M}$ a Riemannian manifold, then the Riemannian risk is exactly the MSE.\hfill$\square$
\end{remark}

\section{Intrinsic Bias of the Estimators}

\subsection{Intrinsic Bias of the Sample Fr\'echet Mean}
\begin{theorem}
Let $\boldsymbol{R}_1,\boldsymbol{R}_2,\ldots,\boldsymbol{R}_N$, be independent and identically distributed (IID)  random matrices with  $K\bR \,{\displaystyle\EqualDist}\,
 {\cal W}_p^C(K,{\bI}),$
and let $\hat{\bP}_{\bI}$ be their sample Fr{\'e}chet mean. Then
\begin{equation}\label{eq:ibiasone}
\text{ibias}(\hat{\bP}_{\boldsymbol{I}})=p\log^2(a),
\end{equation}
where 
\begin{equation}\label{eq:defineaa}
a=\frac{1}{K}\left(\exp\left\{\sum_{j=1}^p\psi(K-p+j)\right\}  \right)^{1/p},
\end{equation}
and $\psi(\cdot)$ denotes the digamma function. 
\end{theorem}
\begin{IEEEproof}
From \eqref{eq:defibias}, when $\bSigma=\bI,$
the intrinsic bias is given by
\begin{equation}\label{eq:wanta}
\text{ibias}(\hat{\bP}_{\bI})\EqualDef \|E
\{ \text{Lmap}_{\boldsymbol{I}}(\hat{\bP}_{\bI})\}\|^2_{\boldsymbol{I}},
\end{equation}
where 
\begin{equation}\label{eq:wantb}
\text{Lmap}_{\bI}(\hat{\bP}_{\bI})=\bI^{1/2}\log(\bI^{-1/2}\hat{\bP}_{\bI}\bI^{-1/2})\bI^{1/2}=\log(\hat{\bP}_{\bI}).
\end{equation}
From \eqref{eq:posdef} $\hat{\bP}_{\bI}$  is positive definite and Hermitian matrix. From \eqref{eq:deflogm},  $\log(\hat{\bP}_{\bI})$ is Hermitian.
Hence,
$E\{\log(\hat{\bP}_{\bI})\}$ is Hermitian. 
Let
\begin{equation}\label{eq:Elog}
E\{\log(\hat{\bP}_{\bI})\}\EqualDef
\begin{bmatrix}
    d_1    & b_{12}  & \cdots & b_{1p}   \\
    b^*_{12} & d_2     & \cdots & b_{2p}   \\
		\vdots & \vdots  & \ddots & \vdots   \\
    b^*_{1p} & b^*_{2p}  & \cdots & d_p 
\end{bmatrix},
\end{equation} 
where $b^*_{ij}$ is the complex conjugate of $b_{ij}$. For fixed $1\leq i< j\leq p$, we define the elementary matrix 
$\boldsymbol{C}_{(ij)}$ as the matrix formed from $\bI$ by exchanging row $i$ with row $j.$ Then, $\boldsymbol{C}_{(ij)}\boldsymbol{C}_{(ij)}^H=\boldsymbol{I}$ and thus $K\boldsymbol{C}_{(ij)}\bR\,\boldsymbol{C}_{(ij)}^H{\displaystyle\EqualDist}\mathcal{W}_p^{C}(K,\bC_{(ij)}\bI\bC_{(ij)}^H)=\mathcal{W}_p^{C}(K,\bI)$. So 
$
\bC_{(ij)}\bR\,\bC_{(ij)}^H{\displaystyle \EqualDist} \bR.
$
 Since $\bC_{(ij)}\in GL(p)$ we know from the equivariant property  that the sample Fr{\'e}chet mean of $\bC_{(ij)}\bR_k\bC_{(ij)}^H, k=1,\ldots,N,$ is  $\bC_{(ij)}\hat{\bP}_\bI\bC_{(ij)}^H$. So,
\begin{equation}\label{eq:equaldist}
\bC_{(ij)}\hat{\bP}_\bI\bC_{(ij)}^H \EqualDist \hat{\bP}_\bI.
\end{equation}
The expression on the left of \eqref{eq:equaldist} is clearly Hermitian, and it is congruent to $ \hat{\bP}_\bI$ and $\bC_{(ij)}$ has full-rank, so $\bC_{(ij)}\hat{\bP}_\bI\bC_{(ij)}^H \in \pdH.$ 
Therefore its logarithm exists as defined in \eqref{eq:deflogm}.
Using \eqref{eq:equaldist} we have 
\begin{eqnarray}
E\{\log(\hat{\bP}_{\boldsymbol{I}})\}&=E\{\log(\boldsymbol{C}_{(ij)}\hat{\bP}_{\boldsymbol{I}}\boldsymbol{C}_{(ij)}^H)\}\label{eq:onea}\\
                                      &=E\{\boldsymbol{C}_{(ij)}\log(\hat{\bP}_{\boldsymbol{I}})\boldsymbol{C}_{(ij)}^H\}\label{eq:oneb}\\
		&=\boldsymbol{C}_{(ij)}E\{\log(\hat{\bP}_{\boldsymbol{I}})\}\boldsymbol{C}_{(ij)}^H.\nonumber
\end{eqnarray}
To make the step from \eqref{eq:onea} to \eqref{eq:oneb} we used 
\eqref{eq:onec} with $\bA$ equated to $\boldsymbol{C}_{(ij)}^H.$
Now $\diag({ E\{\log(\hat{\bP}}_\bI)\})=(d_1\ldots,d_i,\ldots,d_j,\ldots, d_p),$ and $\diag(\bC_{(ij)}E\{\log(\hat{\bP}_\bI)\}\bC_{(ij)}^H)=(d_1\ldots,d_j,\ldots,d_i,\ldots, d_p),$ so we conclude that $d_i=d_j.$ Since this is true for all 
$1\leq i< j\leq p,$ we must have
$d_i=a_0, i=1,\ldots,p,$ for some factor $a_0.$

With reference to \eqref{eq:Elog}, we now demonstrate that $b_{ij} = b^*_{ij} = 0$, $1\leq i < j\leq p.$ For a fixed $i$, we define the elementary matrix,  $\boldsymbol{C}_{(i)},$ which is identical to $\bI$ except  the $i$th row of $\boldsymbol{I}$ is multiplied by $-1$. Note that $\boldsymbol{C}_{(i)}$ is both Hermitian and full-rank, just like $\boldsymbol{C}_{(ij)}.$
So, applying the same reasoning as above,
$
\bC_{(i)}E\{\log(\hat{\bP}_\bI)\}\bC_{(i)}^H
=E\{\log(\hat{\bP}_\bI)\}.
$ 
But the first $i-1$ elements of the $i$th column of $E\{\log(\hat{\bP}_\bI)\}$ are $[b_{1i},\ldots,b_{i-1,i}]^T,$
while the first $i-1$ elements of the $i$th column of 
$\bC_{(i)}E\{\log(\hat{\bP}_\bI)\}\bC_{(i)}^H$ are $[-b_{1i},\ldots,-b_{i-1,i}]^T,$ so that $[b_{1i},\ldots,b_{i-1,i}]^T = {\bf 0}.$ Since this holds for all $i=2,\ldots,p,$  we have that $b_{ij}=0, 1\leq i < j\leq p,$ and their conjugates are likewise zero.
We thus conclude that 
\begin{equation}\label{eq:resulta}
E\{\log(\hat{\bP}_{\boldsymbol{I}})\}=a_0 \bI_p.
\end{equation}

In \cite[p.~4560]{ZhuangWalden17} it was proved that $E\{\log(|\hat{\bP}_{\boldsymbol{I}}|)=p\log(a),$ 
So, denoting the eigenvalues of $\log(\hat{\bP}_{\bI})$ by $\lambda_1,\ldots,\lambda_p$, 
\begin{equation*}
\begin{split}
pa_0 &= \tr( E\{\log(\hat{\bP}_\bI)\} )=E\tr\{ \log(\hat{\bP}_{\boldsymbol{I}})\}\\
		 &=E\{\log\lambda_1+\ldots+\log\lambda_p\}\\
		 &=E\{\log(|\hat{\bP}_{\boldsymbol{I}}|)\}= p\log(a).
\end{split}
\end{equation*}
Hence, $a_0=\log(a),$ which depends on $K$ and $p,$ and
\begin{equation*}
\!\!\!\!\text{ibias}(\hat{\bP}_{\boldsymbol{I}})=\|a_0\bI\|_{\bI}^2=\tr\{a_0\bI a_0\bI\}=pa_0^2=p\log^2(a).
\end{equation*}
\end{IEEEproof}
\begin{remark}
From \eqref{eq:wanta}, \eqref{eq:wantb} and \eqref{eq:wantc} we see that we can write
\begin{equation}\label{eq:wantd}
\text{ibias}(\hat{\bP}_{\boldsymbol{I}})=\|E\{\log(\hat{\bP}_{\boldsymbol{I}})\}\|_{\rm Fr}^2=\tr( E^2\{\log(\hat{\bP}_{\boldsymbol{I}})\}).
\end{equation}
\end{remark}
\begin{remark}\label{remark:noN} 
Note that $\text{ibias}(\hat{\bP}_{\boldsymbol{I}})$ does not depend on the sample size $N$. 
To demonstrate result \eqref{eq:resulta} we carried out some simple simulations 
using $p=3$, $N=3$, $K=20$, with $10\,000$ replications. 
The vector of averaged diagonal entries obtained was
$
[-0.0795, -0.0792, -0.0767]
$
while $a_0=\log(a) = -0.0788.$ 
Repeating with  the same parameters, except $N=10$, we obtained instead
$
[-0.0794, -0.0783, -0.0796].
$
Off-diagonal terms were almost zero in both cases.
\hfill$\square$
\end{remark}

We now turn to the more general case.
Let $\bS_1,\ldots,\bS_N$ be IID samples with $K\boldsymbol{S}\,{\displaystyle\EqualDist}\,\mathcal{W}_p^{C}(K,\bSigma),$
and consider their sample Fr{\'e}chet mean $\hat{\bP}_{\bSigma}.$ Then $\text{ibias}(\hat{\bP}_{\boldsymbol{\Sigma}})$ is given by
\begin{eqnarray}
&&\!\!\!\!\!\!\!\!\!\!\!\!\!\!\!\!\!\!\!\!\!
\|E\{ \text{Lmap}_{\boldsymbol{\Sigma}}(\hat{\bP}_{\boldsymbol{\Sigma}})\}\|^2_{\boldsymbol{\Sigma}}\nonumber\\
=&&\!\!\!\!\!\!\!\!\|E\{ \boldsymbol{\Sigma}^{1/2}\log(\boldsymbol{\Sigma}^{-1/2}\hat{\bP}_{\boldsymbol{\Sigma}}\boldsymbol{\Sigma}^{-1/2})\boldsymbol{\Sigma}^{1/2}\}\|^2_{\boldsymbol{\Sigma}}\nonumber\\
=&&\!\!\!\!\!\!\!\!\|\boldsymbol{\Sigma}^{1/2}E\{\log(\boldsymbol{\Sigma}^{-1/2}\boldsymbol{\Sigma}^{1/2}\hat{\bP}_{\bI}\boldsymbol{\Sigma}^{1/2}\boldsymbol{\Sigma}^{-1/2})\}\boldsymbol{\Sigma}^{1/2}\|^2_{\boldsymbol{\Sigma}}\nonumber\\
=&&\!\!\!\!\!\!\!\!\|\boldsymbol{\Sigma}^{\half}E\{\log(\hat{\bP}_{\bI})\}\boldsymbol{\Sigma}^{\half}\|^2_{\boldsymbol{\Sigma}}
\nonumber\\
=&&\!\!\!\!\!\!\!\!\|a_0\boldsymbol{\Sigma}\|^2_{\boldsymbol{\Sigma}}=\text{tr}(\bSigma^{-1}a_0\bSigma\bSigma^{-1}a_0\bSigma)=pa_0^2(K).\label{eq:final}
\end{eqnarray}
So the intrinsic bias of the sample Fr\'echet mean of the complex Wishart matrices with a general true covariance matrix $\boldsymbol{\Sigma}$ does not depend on the parameter $\boldsymbol{\Sigma},$ nor on $N,$ but just on $p$ and $K.$ We particularly note that
\begin{equation}\label{eq:note}
E\{\text{Lmap}_\bSigma(\hat{\bP}_{\boldsymbol{\Sigma}})\}=a_0\boldsymbol{\Sigma}.
\end{equation}
\subsubsection{Single Covariance Matrix}
As a special case, by taking $N=1$ in \eqref{eq:note}, we see that
\begin{equation}\label{eq:smith}
E\{\text{Lmap}_\bSigma({\bS})\}=a_0\boldsymbol{\Sigma}.
\end{equation}
Using the result $\psi (z+1)=\psi (z)+(1/z)$ it can be shown after some manipulation that 
 the form of \eqref{eq:smith}
matches that in \cite[eqns.~(100), (102)]{Smith05}.

\subsubsection{Scalar Case}
Consider the scalar case, $(p=1),$
where $KS\,{\displaystyle\EqualDist}\,\mathcal{W}_1^{C}(K,1),$ i.e., the single entry of the covariance matrix is unity. From \eqref{eq:final} we have that
$\text{ibias}(S)=a_0^2$  so that, using
\eqref{eq:defineaa},
\begin{equation}\label{eq:scalarchi}
\text{ibias}(S)=\log^2(a)=\left[ \psi(K)-\log(K)\right]^2.
\end{equation}
It is straightforward to show that $S\,{\displaystyle\EqualDist}\,\frac{1}{2K}\chi^2_{2K},$ (scaled chi-square with $2K$ degrees of freedom), so that \eqref{eq:scalarchi} gives the intrinsic bias in this case. 
Now suppose $Y_1,\ldots, Y_K \,{\displaystyle\EqualDist}\,{\rm Exp}(1)\,{\displaystyle\EqualDist}\,\frac{1}{2}\chi^2_{2}.$ Then ${\overline Y},$ the mean, has distribution ${\overline Y}\,{\displaystyle\EqualDist}\,\frac{1}{2K}\chi^2_{2K},$ and so
$\text{ibias}({\overline Y})$ is given by \eqref{eq:scalarchi}; this result agrees with 
\cite[p.~1569]{OllerCorcuera95}.

\subsection{Intrinsic Bias of the Sample Mean}\label{sec:average}
Let $\bS_1,\ldots,\bS_N$ be IID samples with $K\boldsymbol{S}\,{\displaystyle\EqualDist}\,\mathcal{W}_p^{C}(K,\bSigma),$
and consider their ordinary sample  mean $\bar{\bS}.$ 
We know that
\begin{equation}\label{eq:composite}
KN {\bar{\bS}}
 \EqualDist
{\cal W}_p^C(KN,{\bSigma}).
\end{equation}
The intrinsic bias of the sample mean can be inherited from that of the Fr{\'e}chet mean by taking the sample size of the Fr{\'e}chet mean
to be $N=1,$ and then replacing the degrees of freedom $K$ by $KN$ in \eqref{eq:final}:
\begin{equation}\label{eq:simpleav}
\text{ibias}( {\bar{\bS}} )=
\|E\{ \text{Lmap}_{\boldsymbol{\Sigma}}({\bar{\bS}})\}\|^2_{\boldsymbol{\Sigma}}=
p \,a_0^2(KN),
\end{equation}

\section{Intrinsic Risk of the Estimators}
\subsection{Intrinsic Risk of the Sample Fr\'echet Mean}
The Riemannian or intrinsic risk was defined in 
\eqref{eq:defRisk}. We now examine this for 
the sample Fr\'echet mean.
\begin{eqnarray*}
&&E\{ \|\text{Lmap}_{\bSigma}(\hat{\bP}_{\bSigma})\|^2_{\bSigma}\}\\
&&=E\{ \|\boldsymbol{\Sigma}^{1/2}\log(\boldsymbol{\Sigma}^{-1/2}\hat{\bP}_{\boldsymbol{\Sigma}}\boldsymbol{\Sigma}^{-1/2})\boldsymbol{\Sigma}^{1/2}\|^2_{\bSigma}\}\\
&&=E\{\| \boldsymbol{\Sigma}^{1/2}\log(\boldsymbol{\Sigma}^{-1/2}\boldsymbol{\Sigma}^{1/2}\hat{\bP}_{\bI}\boldsymbol{\Sigma}^{1/2}\boldsymbol{\Sigma}^{-1/2})\boldsymbol{\Sigma}^{1/2}
\|^2_{\bSigma}\}\\
&&=E\{\| \boldsymbol{\Sigma}^{1/2}\log( \hat{\bP}_{\bI})\boldsymbol{\Sigma}^{1/2}\|^2_{\bSigma}\}\\
&&=E\{\tr(\log^2( \hat{\bP}_{\bI}) ) \},
\end{eqnarray*}
where for the last step we have used \eqref{eq:scFrob} and the cyclic nature of trace. From \eqref{eq:cdista} and \eqref{eq:cdistb}  we alternatively can write
\begin{eqnarray*}
E\{ \|\text{Lmap}_{\bSigma}(\hat{\bP}_{\bSigma})\|^2_{\bSigma}\}\!\!\!\!&=&\!\!\!\!E\{ d_g^2(\bI,\hat{\bP}_{\bI})\}=E\|\log(\hat{\bP}_{\bI})\|^2_{\rm Fr}.
\end{eqnarray*}

Hence the Riemannian risk is independent of the underlying true covariance matrix. The Riemannian risk can be decomposed into two parts, (i) the intrinsic bias which does not depend on the sample size, and (ii) the sum of the variances of every  entry of  $\text{Lmap}_{\bI}(\hat{\bP}_{\bI})$. Let $L_{ij}=(\log(\hat{\bP}_{\bI}))_{ij},$ the $(i,j)$th entry of the Hermitian matrix $\log(\hat{\bP}_{\bI}).$ Then
\begin{equation*}
\begin{split}
E\{ \|\text{Lmap}_{\bSigma}(\hat{\bP}_{\bSigma})\|^2_{\bSigma}\}&=E\|\log(\hat{\bP}_{\bI})\|^2_{\rm Fr}= \sum_{i,j=1}^p\ E\{ |L_{ij}|^2\}\\
&= \sum_{i,j=1}^p\var\{L_{ij}\}+ \sum_{i,j=1}^p |E\{ L_{ij}\}|^2\\
&= \sum_{i,j=1}^p\var\{L_{ij}\}+ \|E\{\log(\hat{\bP}_{\boldsymbol{I}})\}\|_{\rm Fr}^2\\
&=\!\sum_{i,j=1}^p\!\!\var\{(\log(\hat{\bP}_{\bI}))_{ij}\}+\text{ibias}(\hat{\bP}_{\bI}),
\end{split}
\end{equation*} 
where we have used \eqref{eq:wantd} for the last line. 
\subsection{Intrinsic Risk of the Sample Mean}
If we set $N=1,$ and replace the degrees of freedom $K$ by $KN,$ we can replace $\hat{\bP}_{\bSigma}$ by $\bar{\bS}$ and $\hat{\bP}_{\bI}$ by $\bar{\bS}_\bI$ to give
$$
E\{ \|\text{Lmap}_{\bSigma}(\bar{\bS})\|^2_{\bSigma}\}=\sum_{i,j=1}^p\var\{(\log(\bar{\bS}_{\bI}))_{ij}\}+\text{ibias}(\bar{\bS}_\bI).
$$

\section{Comparisons}
\subsection{Intrinsic Bias}
We firstly examine $a_0(K)$ for $K\geq p.$ From \eqref{eq:defineaa} we have 
$
a_0(K)=(1/p)\sum_{j=1}^p \psi(K+1-j)-\log K.
$
We know $\psi(K)-\log K<0$ and $\psi(K)$ increases with positive $K,$ so $\psi(K-1)-\log K<0,\ldots,\psi(K+1-p)-\log K<0.$ Therefore,
$
a_0(K)=(1/p)[\sum_{j=1}^p \psi(K+1-j)-p\log K]<0.
$
Using $\psi (z+1)=\psi (z)+(1/z)$ we have
\begin{equation*}
\begin{split}
a_0(K+1)&=\frac{1}{p}\sum_{j=1}^p\left[\psi(K+1-j)+\textstyle{\frac{1}{K+1-j}}\right]-\log(K+1)\\
&=a_0(K) + \frac{1}{p}\sum_{j=1}^p\left[\textstyle{\frac{1}{K+1-j}}\right]
-\log\left(1+{\textstyle{\frac{1}{K}}}\right).
\end{split}
\end{equation*}
Clearly $(1/p)\sum_{j=1}^p\left[\textstyle{\frac{1}{K+1-j}}\right]>(1/K)$ and $(1/K)-\log\left(1+{\textstyle{\frac{1}{K}}}\right) > 0$ for $K\geq 1$ \cite[p.~159]{Bullen98}, so $a_0(K+1) > a_0(K),$ and since $a_0(K)$ is negative, it is monotonically increasing and bounded above by zero, and $a_0^2(K)$ decreases in magnitude with increasing $K.$
So, for $N\geq 2,$
\begin{equation}\label{eq:avbetter}
\text{ibias}(\hat{\bP}_{\boldsymbol{\Sigma}})=
\text{ibias}(\hat{\bP}_{\boldsymbol{I}})=pa_0^2(K) >pa_0^2(KN)=\text{ibias}( {\bar{\bS}} ),
\end{equation}
i.e., the simple average is preferable in terms of intrinsic bias for $N\geq 2.$
\subsubsection{Scalar case}
Set $p=1$ to obtain
for $N\geq 2,$
\begin{equation}\label{eq:avbetterscal}
\text{ibias}(S)=a_0^2(K) >a_0^2(KN)=\text{ibias}( {\bar{S}} ),
\end{equation}
Since $\text{ibias}( {\bar{S}} )=\left[ \psi(KN)-\log(KN)\right]^2,$ (see \eqref{eq:scalarchi}), and $\psi(z)\rightarrow \log z$ as $z\rightarrow\infty,$  we see that $\lim_{N\rightarrow\infty} \text{ibias}( {\bar{S}} )=0.$
\subsection{Riemannian Risk}

We know that the sample Fr{\'e}chet mean $\hat{\bP}_{\bI}$ converges almost surely to the population 
Fr{\'e}chet mean ${\bP}_{\bI}$ \cite[p.~4555]{ZhuangWalden17}. So we can conclude that
$$
\lim_{N\rightarrow \infty} E\{ \|\text{Lmap}_{\bSigma}(\hat{\bP}_{\bSigma})\|^2_{\bSigma}\}=\text{ibias}(\hat{\bP}_{\bI})=p\a_0^2(K).
$$
For the sample mean we know $\bar{\bS}_{\bI}$ converges almost surely to the population mean $\bI$. So we have 
\begin{equation*}
\begin{split}
\lim_{N\rightarrow \infty}E\{ \|\text{Lmap}_{\bSigma}(\bar{\bS})\|^2_{\bSigma}\}&=\lim_{N\rightarrow \infty}\text{ibias}(\bar{\bS}_\bI)\\
&=\lim_{N\rightarrow \infty}pa_0^2(KN).
\end{split}
\end{equation*}

From \eqref{eq:avbetter} we thus see that, {\it asymptotically}, the Riemannian risk will be smaller for the simple average.
\subsubsection{Scalar case}
For $p=1$ and finite $N$, the Riemannian risk is $\var\{\log(\hat{P}_{1})\}+E^2\{\log(\hat{P}_1)\},
$  with
$
\log({\hat P}_1)=(1/N)\sum_{j=1}^N \log(S_j).
$
Then (e.g., \cite[p.~275]{PercivalWalden00}), 
$E\{ \log(S_j)\}=\psi(K)-\log(K)$ and $\var\{ \log(S_j)\}=\psi'(K),$ where $\psi'(\cdot)$ is the trigamma function.
So the Riemannian risk of
${\hat P}_1,$ ${\rm Rr}({\hat P}_1),$ say,  is given by
\begin{equation}\label{eq:theorisk}
{\rm Rr}({\hat P}_1)=\frac{1}{N} \psi'(K) +\left[ \psi(K)-\log(K)\right]^2.
\end{equation}
When $N=1,$ \eqref{eq:theorisk} gives
$
{\rm Rr}(S)=\psi'(K) +\left[ \psi(K)-\log(K)\right]^2
$
so, replacing $K$ by $KN,$
\begin{equation}\label{eq:riskmn}
{\rm Rr}({\bar S})=\psi'(KN) +\left[ \psi(KN)-\log(KN)\right]^2.
\end{equation}
Plotting the difference ${\rm Rr}({\hat P}_1)-{\rm Rr}({\bar S})$  as a function of $N$ and $K,$ we find that  the risk is always higher for ${\rm Rr}({\hat P}_1).$
\section*{Acknowledgement}
Linjie Zhuang was supported by an Imperial College PhD Scholarship.

\section{Conclusions}
We have proved that 
the simple average of $N$ Wishart matrices is preferred over the Fr\'echet mean (i)
in terms of intrinsic bias for $N\geq 2,$ and (ii) likewise for the asymptotic Riemannian risk, and for the  Riemannian risk itself in the scalar case. Simple expressions were given for the intrinsic bias for both types of average. Non-intrinsic performance measures in \cite{ZhuangWalden17} also favoured the simple average. There is thus strong evidence to prefer the sample  mean in this context.


\end{document}